\documentclass[12pt]{article}
\usepackage{mathrsfs}
\usepackage{graphicx,amsmath,amsfonts,amssymb,amscd,amsthm}
\newtheoremstyle{kai}
{3pt} {3pt} {} {} {\bfseries} {.} {.5em} {}
\def\EquationsBySection{\def\theequation
{\thesection.\arabic{equation}}%
\@addtoreset{equation}{section}}

\newcommand\old[1]{}
\newcommand{\pend}{\hfill \thicklines \framebox(6.6,6.6)[l]{}}
\renewenvironment{proof}{\noindent {\it  Proof.} \rm}{\pend}
\newtheorem{theorem}{Theorem}[section]
\newtheorem{lemma}{Lemma}[section]

\newtheorem{remark}{Remark}[section]

\EquationsBySection \makeatother
\usepackage{indentfirst}
\begin{document}
\pagestyle{plain}
\title
{\bf Rate of convergence for numerical solutions to SFDEs with
jumps}
\author{
Jianhai Bao$^1$, Xuerong Mao$^2$, Chenggui Yuan$^{3}$\thanks{{\it E-mail address: }    jianhaibao@yahoo.com.cn, xuerong@stams.strath.ac.uk,
C.Yuan@swansea.ac.uk.}
\\
\\
$^1$School of Mathematics, Central South University,\\
 Changsha, Hunan
410075, P.R.China
 \\
  $^2$Department of Statistics and Modelling Science,\\
University of Strathclyde, Glasgow G1 1XH, UK
\\
  $^3$Department of Mathematics, Swansea University,\\
 Swansea SA2 8PP, UK
 }

\date{}
\maketitle
\begin{abstract}{\rm In this paper, we are interested in the
numerical solutions of stochastic functional differential equations (SFDEs) with
{\it jumps}.  Under the global Lipschitz condition, we show that
the $p$th moment convergence of the Euler-Maruyama (EM)
numerical solutions to SFDEs with jumps has order $1/p$ for any $p\ge 2$.
This is significantly different from the case of SFDEs without jumps
where the order is $1/2$ for any $p\ge 2$.  It is therefore best
to use the mean-square convergence for SFDEs with jumps.
Consequently, under the local Lipschitz condition,
 we reveal that the order of the mean-square convergence is close to $1/2$,
provided that the local
Lipschitz constants, valid on balls of radius $j$, do not
grow faster than $\log j$.
}\\

\end{abstract}
\noindent

\section{Introduction}
Recently, the theory of functional differential equations (FDEs) has
received a great deal of attention. Hale and Lune \cite{Hal} have
studied deterministic functional differential equations (DFDEs) and
their stability. For stochastic functional differential equations
(SFDEs), we here highlight the great contribution of Kolmanovskii
and Nosov \cite{Kol} and Mao \cite{m97}. Kolmanovskii and Nosov
\cite{Kol} not only established the theory of existence and
uniqueness of SFDEs but also investigated the stability and
asymptotic stability of the equations, while Mao \cite{m97} studied
the exponential stability of the equations.

On the other hand, stochastic differential equations (SDEs) with
jumps have been widely used in many branches of science and industry,
in particular, in economics,
finance and engineering (see, for
example, Gukhal \cite{g04}, R. Cont \cite{ct04}, Sobczyk \cite{s91}
and references therein).   Since most SDEs with jumps cannot  be solved
explicitly,  numerical methods have become essential.
Under the local Lipschitz condition, Higham and Kloeden
\cite{hk05} showed the strong convergence and nonlinear stability
for the EM numerical solutions to SDEs with jumps, while, in
\cite{hk07}, Higham and Kloeden further revealed the strong
convergence rate for the backward Euler (BE) on SDEs with jumps,
provided that the drift coefficient obeys one-side Lipschitz
condition and polynomial growth condition.

Returning to the SFDEs, we recalled
Mao \cite{m03} developed a numerical scheme for them.
 Under the {\it local
Lipschitz condition}, Mao \cite{m03} showed the strong convergence
of the EM numerical solutions to SFDEs, but revealed the rate of
the convergence  under the {\it global Lipschitz condition}.
But there is so far no work on numerical methods for SFDEs
with jumps.

Motivated by the papers mentioned above, we are here interested in the
numerical solutions to SFDEs {\it
with jumps}.  In comparison with the results obtained by
Mao \cite{m03}, our significant contributions are:

\begin{itemize}

\item  Under the global Lipschitz condition, we show that
the $p$th moment convergence of the EM
numerical solutions to SFDEs with jumps has order $1/p$ for any $p\ge 2$.
This is significantly different from the case of SFDEs without jumps
where the order is $1/2$ for any $p\ge 2$.  In practice, it is therefore best
to use the mean-square convergence for SFDEs with jumps.

\item Under the local Lipschitz condition, Mao \cite{m03}
showed the strong convergence \emph{without rate}
of the EM numerical solutions to SFDEs without jumps.
However, we shall reveal that the order of the mean-square convergence is closed to $1/2$,
provided that the local
Lipschitz constants, valid on balls of radius $j$, do not
grow faster than $\log j$.  More precisely, the order of the mean-square convergence is  $1/(2+\epsilon)$,
provided that the local
Lipschitz constants do not
grow faster than $(\log j)^{1/(1+\epsilon)}$.

\item  Some new techniques are developed to cope with
the difficulty due to the jumps.

\end{itemize}

This paper is organized as follows: Section 2 gives
some preliminary results, in particular, the EM numerical solutions
to SFDEs with jumps are set up. In section 3,
we discuss the $p$th moment convergence of the EM
numerical solutions to SFDEs with jumps under the global Lipschitz
condition.
The rate of the mean-square convergence of the EM
numerical solutions to SFDEs with jumps under the local Lipschitz
condition is provided in Section 4.  Finally, in order to make the
paper self-contained, an existence-and-uniqueness result of
solutions to SFDEs with jumps is provided in the Appendix.

\section{Preliminaries}
Throughout this paper, we let $\{\Omega,{\mathcal F},\{{\mathcal
F}_{t}\}_{t\geq0}, P\}$ be a complete probability space with a
filtration $\{{\mathcal F}_{t}\}_{t\geq0}$ satisfying the usual
conditions (i.e., it is continuous on the right and
$\mathcal{F}_{0}$-contains all $P$-zero sets). Let $|\cdot|$ denote
the Euclidean norm and the matrix trace norm. Let $\tau>0$ and
$D:=D([-\tau,0];R^n)$ denote the family of all right-continuous
functions with left-hand limits $\varphi$ from $[-\tau,0]$ to $R^n,$ and
$\hat D:=\hat D([-\tau,0];R^n)$ denote the family of all left-continuous
functions with right-hand limits $\varphi$ from $[-\tau,0]$ to $R^n,$
we will always use
$\|\varphi\|=\sup_{-\tau\leq\theta\leq0}|\varphi(\theta)|$ to denote the norm in $D$ and $\hat D$ potentially involved when no confusion possibly arises.
$D^b_{\mathcal {F}_0}([-\tau,0];R^n)$ denotes the family of all
almost surely bounded, ${\mathcal {F}}_0$-measurable,
$D([-\tau,0];R^n)$-valued random variables. For all $t\geq0$,
$x_t=\{x(t +\theta):-\tau\leq\theta\leq0\}$ is regarded as a
$D([-\tau,0];R^n)$-valued stochastic process. Let $x(t^-)$
denotes $\lim_{s\uparrow t}x(s)$ on $t\in[-\tau, T]$ and $x_{t^-}=\{x(t+\theta)^-:-\tau\leq\theta\leq0\},$ it is easy to see $x(t^-)$ is a $\hat D([-\tau,0];R^n)$-valued stochastic process.

In this paper, we consider the following SFDE with jumps
\begin{equation}\label{eq1}
dx(t)=f(x_{t^-})dt+g(x_{t^-})dB(t)+h(x_{t^-})dN(t) ,\;\;\;\;\; 0\le
t\le T,
\end{equation}
with the initial data $x_{0}=\xi\in D^{b}_{{\mathcal
F}_{0}}([-\tau,0];R^n)$.  Here, $f, h:
\hat D([-\tau,0];R^n)\rightarrow R^n$, $g:\hat D([-\tau,0];R^n)\rightarrow
R^{n\times m}$, $B(t)$ is an $m$-dimensional Brownian motion and
$N(t)$ is a scalar Poisson process with intensity $\lambda.$ We
further assume that   $B(t)$ and $N(t)$ are independent. It should be pointed out that the solution of Eq. \eqref{eq1} is in $D([-\tau,0];R^n)$

For our purposes, we need the following assumptions which can also
guarantee the existence and uniqueness of solution to \eqref{eq1}
(see Appendix).

\begin{enumerate}
\item[\textmd{(\bf H1)}] (Global Lipschitz condition) There exists a left-continuous nondecreasing function
$\mu:[-\tau,0]\rightarrow R_+$ such that, for all $\varphi,\psi \in
\hat D([-\tau,0];R^n)$,
\begin{equation}\label{eq2}
|f(\varphi)-f(\psi)|^2\vee|g(\varphi)-g(\psi)|^2\vee|h(\varphi)-h(\psi)|^2\leq
\int_{-\tau}^0|\varphi(\theta)-\psi(\theta)|^2d\mu(\theta).
\end{equation}
\end{enumerate}

\begin{remark}
For simplicity, we write $L=\mu(0)-\mu(-\tau)$, which is referred to
as the global Lipschitz constant. Note from \eqref{eq2} that, for
all $\varphi,\psi \in \hat D([-\tau,0];R^n)$,
\begin{equation}\label{eq4}
|f(\varphi)-f(\psi)|^2\vee|g(\varphi)-g(\psi)|^2\vee|h(\varphi)-h(\psi)|^2\leq
L\|\varphi-\psi\|^2.
\end{equation}
This further implies the linear growth condition; that is, for
$\varphi\in \hat D([-\tau,0];R^n)$,
\begin{equation}\label{eq5}
|f(\varphi)|^2\vee |g(\varphi)|^2\vee |h(\varphi)|^2\leq
K(1+\|\varphi\|^2),
\end{equation}
where $K=2(L\vee|f(0)|^2\vee|g(0)|^2\vee|h(0)|^2)$.

\end{remark}

\begin{enumerate}
\item[\textmd{(\bf H2)}] (Continuity of initial data) For $\xi\in D^{b}_{{\mathcal
F}_{0}}([-\tau,0];R^n), 0 \le u\le \tau$ and $p\geq2$, there is a constant $\beta>0$
such that
\begin{equation}\label{eq6}
E\left(\sup\limits_{-\tau\leq s\leq t\leq
0, \atop
|t-s|\le u}|\xi(s)-\xi(t)|^p\right)\leq\beta u.
\end{equation}
\end{enumerate}

For given $T\geq0$ and $\tau>0$, the time-step size
$\bigtriangleup\in (0,1)$ is defined by
\begin{equation*}
\bigtriangleup=\frac{\tau}{N}=\frac{T}{M}
\end{equation*}
with some integers $N>\tau$ and $M>T$. The EM method applied to
\eqref{eq1} produces approximations $\bar{y}(k\bigtriangleup)\approx
x(k\bigtriangleup)$  by setting
$\bar{y}(k\bigtriangleup)=\xi(k\bigtriangleup),-N\leq k\leq 0,$ and
\begin{equation}\label{eq7}
\bar{y}((k+1)\bigtriangleup)=\bar{y}(k\bigtriangleup)
+f(\bar{y}_{k\bigtriangleup})\bigtriangleup+g(\bar{y}_{k\bigtriangleup})\bigtriangleup
B_k+h(\bar{y}_{k\bigtriangleup})\bigtriangleup N_k,
\end{equation}
where $\bigtriangleup B_k=B((k+1)\bigtriangleup)-B(k\bigtriangleup)$
is a Brownian increment, $\bigtriangleup
N_k=N((k+1)\bigtriangleup)-N(k\bigtriangleup)$ is a Poisson
increment,  and
$\bar{y}_{k\bigtriangleup}=\{\bar{y}_{k\bigtriangleup}(\theta):-\tau\leq\theta\leq0\}$
is a $D([-\tau,0];R^n)$-valued random variable defined by
\begin{equation}\label{eq8}
\bar{y}_{k\bigtriangleup}(\theta)=\frac{(i+1)\bigtriangleup-\theta}{\bigtriangleup}\bar{y}((k+i)\bigtriangleup)+
\frac{\theta-i\bigtriangleup}{\bigtriangleup}\bar{y}((k+i+1)\bigtriangleup)
\end{equation}
for $i\bigtriangleup\leq\theta\leq(i+1)\bigtriangleup$, $i=-N,
-(N-1), \cdots, -1$, where in order for $\bar{y}_{-\bigtriangleup}$
to be well defined, we set
$\bar{y}(-(N+1)\bigtriangleup)=\xi(-N\bigtriangleup)$.

Given the discrete-time approximation
$\{\bar{y}(k\bigtriangleup)\}_{k\geq0}$, we define a continuous-time
approximation $y(t)$ by $y(t)=\xi(t)$ for $-\tau\leq t\leq 0$, while
for $t\in[0,T]$,
\begin{equation}\label{eq9}
\begin{split}
y(t)&=\xi(0)+\int_0^tf(\bar{y}_{s^-})ds+\int_0^tg(\bar{y}_{s^-})dB(s)
+\int_0^th(\bar{y}_{s^-})dN(s),
\end{split}
\end{equation}
where, for fixed $\theta\in[-\tau,0]$,
\begin{equation*}
\bar{y}_{t^-}=\lim_{s\uparrow t}\bar{y}_s,\ \ \ \ \ \ \
\bar{y}_t=\sum\limits_{k=0}^{M-1}\bar{y}_{k\bigtriangleup}I_{[k\bigtriangleup,
(k+1)\bigtriangleup)}(t).
\end{equation*}
It is easy to see $y(k\bigtriangleup)=\bar{y}(k\bigtriangleup)$ for
$k=-N, -N+1,\cdots, M$. That is, the discrete-time and
continuous-time EM numerical solutions coincide at the gridpoints.

\begin{remark}
It is easy to observe from \eqref{eq8} that
\begin{equation}\label{eq12}
\| \bar{y}_{k\bigtriangleup}\|=\max\limits_{-N\leq i\leq0}|
\bar{y}((k+i)\bigtriangleup)|, \ \ \ k=-1, 0, 1, \cdots, M-1,
\end{equation}
which further yields
\begin{equation*}
\| \bar{y}_{k\bigtriangleup}\|\leq\| y_{k\bigtriangleup}\|, \ \ \
k=-1, 0, 1, \cdots, M-1,
\end{equation*}
by  $y(k\bigtriangleup)=\bar{y}(k\bigtriangleup)$ and, for any
$t\in[0,T]$,
\begin{equation}\label{eq14}
\|
\bar{y}_t\|=\|\bar{y}_{[\frac{t}{\bigtriangleup}]\bigtriangleup}\|
\leq\|y_{[\frac{t}{\bigtriangleup}]\bigtriangleup}\|\leq\sup\limits_{-\tau\leq
s\leq t}|y(s)|,
\end{equation}
where $[\frac{t}{\bigtriangleup}]$ is the integer part of
$\frac{t}{\bigtriangleup}$.

\end{remark}

\section{Convergence under the global Lipschitz condition}

In this section, we will investigate the rate of the convergence under the global Lipschitz condition.  Our results
 reveal  a significant difference  from these on the SDEs without jumps.

\begin{lemma}\label{le1}
Under the condition (\ref{eq5}), for any $p\geq2$
there exists a positive constant $H(p):=H(p,T,\xi, K)$ which may dependent on $ p, T, \xi,
K$ such that
\begin{equation}\label{eq10}
E\Big(\sup\limits_{-\tau\leq t\leq T}| x(t)|^p\Big)\vee E\left(\sup\limits_{-\tau\leq t\leq T}| y(t)|^p\right)\leq H(p).
\end{equation}
\end{lemma}

\begin{proof}
Since the arguments of the moment bounds for the exact and continuous
approximate solutions to \eqref{eq1} are very similar, we here only
give an estimate for the continuous approximate solution $y(t)$. For every integer $R\ge 1,$ define the stopping time
$$
\theta_R=\inf\{t\ge 0: \|y_t\|\ge R\}
$$
It is easy to see from \eqref{eq9} that, for any
$t\in[0,T]$,
\begin{equation}\label{eq15}
\begin{split}
&E\left(\sup\limits_{0\leq s\leq t}|
y(s\wedge \theta_R)^-|^p\right)\le E\left(\sup\limits_{0\leq s\leq t}|
y(s\wedge \theta_R)|^p\right)\\
&\leq4^{p-1}\Big[E\|\xi\|^p+E\left(\sup\limits_{0\leq
s\leq t}\left|\int_0^{s}f(\bar{y}_{(r\wedge \theta_R)^-})dr\right|^p\right)\\
&+E\left(\sup\limits_{0\leq s\leq
t}\left|\int_0^{s}g(\bar{y}_{(r\wedge \theta_R)^-})dB(r)\right|^p\right)+E\left(\sup\limits_{0\leq
s\leq t }\left|\int_0^{s}h(\bar{y}_{(r\wedge \theta_R)^-})dN(r)\right|^p\right)\Big].
\end{split}
\end{equation}
Noting that $E\|y_{(t\wedge \theta_R)^-}\|\le R$ and \eqref{eq14}, one may have
 $E\|\bar y_{(t\wedge \theta_R)^-}\|\le R.$
 By the H\"older inequality and \eqref{eq5},
\begin{equation*}
\begin{split}
E\left(\sup\limits_{0\leq s\leq
t}\left|\int_0^sf(\bar{y}_{(r\wedge \theta_R)^-})dr\right|^p\right)&\leq
T^{p-1}\int_0^tE| f(\bar{y}_{(r\wedge \theta_R)^-})|^pdr\\
&\leq T^{p-1}\int_0^tE[K(1+ \|\bar{y}_{(r\wedge \theta_R)^-}\|^2)]^{\frac{p}{2}}dr\\
&=T^{p-1}\int_0^tE[K(1+ \|\bar{y}_{(r\wedge \theta_R)^-}\|^2)]^{\frac{p}{2}}dr\\
&\leq2^{\frac{p}{2}-1}T^pK^{\frac{p}{2}}+2^{\frac{p}{2}-1}T^{p-1}K^{\frac{p}{2}}\int_0^tE\|\bar{y}_{(r\wedge \theta_R)^-}\|^pdr.
\end{split}
\end{equation*}
This, together with \eqref{eq14}, immediately reveals that
\begin{equation}\label{eq16}
E\left(\sup\limits_{0\leq s\leq
t}\left|\int_0^sf(\bar{y}_{(r\wedge \theta_R)^-})dr\right|^p\right)\leq
c_1T+c_1\int_0^tE\left(\sup\limits_{-\tau\leq r\leq
s}|y(r\wedge \theta_R)^-|^p\right)ds,
\end{equation}
where $c_1=2^{\frac{p}{2}-1}T^{p-1}K^{\frac{p}{2}}$. Now, using the
Burkholder-Davis-Gundy inequality \cite[Theorem 7.3, p40]{m97} and
the H\"older inequality, we deduce that there exists a positive
constant $c_p$ such that
\begin{equation*}
\begin{split}
E\left(\sup\limits_{0\leq s\leq
t}\left|\int_0^sg(\bar{y}_{(r\wedge \theta_R)^-})dB(r)\right|^p\right)&\leq c_pE\left(\int_0^t| g(\bar{y}_{(r\wedge \theta_R)^-})|^2dr\right)^{p/2}\\
&\leq c_pT^{\frac{p-2}{2}}\int_0^tE|g(\bar{y}_{r\wedge \theta_R})^-|^pdr.
\end{split}
\end{equation*}
In the same way as \eqref{eq16} was done, it then follows easily
that
\begin{equation*}
\begin{split}
E\left(\sup\limits_{0\leq s\leq
t}\left|\int_0^sg(\bar{y}_{(r\wedge \theta_R)^-})dB(r)\right|^p\right)\leq
c_2T+c_2\int_0^tE\left(\sup\limits_{-\tau\leq r\leq
s}|y(r\wedge \theta_R)^-|^p\right)ds,
\end{split}
\end{equation*}
where $c_2=2^{\frac{p}{2}-1}T^{\frac{p-2}{2}}K^{\frac{p}{2}}c_p$.
Moreover, observing that $\tilde{N}(t)=N(t)-\lambda t, t\geq0$ is a
martingale measure,
using the Burkholder-Davis-Gundy inequality \cite[Theorem 48,
p193]{p04}, H\"older inequality and  \eqref{eq5}, we obtain for some positive constant $\bar{c}_p$,
\begin{equation*}
\begin{split}
&E\left(\sup\limits_{0\leq s\leq
t}\left|\int_0^sh(\bar{y}_{(r\wedge \theta_R)^-})dN(r)\right|^p\right)\\
&\leq
E\left(\sup\limits_{0\leq s\leq
t}\left|\int_0^sh(\bar{y}_{(r\wedge \theta_R)^-})d\tilde{N}(r)+\lambda\int_0^sh(\bar{y}_{(r\wedge \theta_R)^-})dr\right|^p\right)\\
&\leq2^p\left[E\left(\sup\limits_{0\leq s\leq
t}\left|\int_0^sh(\bar{y}_{(r\wedge \theta_R)^-})d\tilde{N}(r)\right|^p+\lambda^p\sup\limits_{0\leq
s\leq t}\left|\int_0^sh(\bar{y}_{(r\wedge \theta_R)^-})dr\right|^p\right)\right]\\
&\leq2^p\Big[\bar{c}_p\lambda^{p/2}E\left(\int_0^t|
h(\bar{y}_{(r\wedge \theta_R)^-}))|^2dr\right)^{p/2}+\lambda^pT^{p-1}\int_0^tE|h(\bar{y}_{(r\wedge \theta_R)^-}))|^pdr\Big]\\
&\leq c_3T+c_3\int_0^tE\left(\sup\limits_{-\tau\leq r\leq
s}|y(r\wedge \theta_R)^-|^p\right)ds,
\end{split}
\end{equation*}
where
$c_3=2^{\frac{3p}{2}-1}K^{\frac{p}{2}}\Big[\bar{c}_p\lambda^{p/2}T^{\frac{p-2}{2}}+\lambda^pT^{p-1}\Big]$.
Hence, in \eqref{eq15}
\begin{align*}
&E\left(\sup\limits_{0\leq s\leq t}|
y(s\wedge \theta_R)^-|^p\right)\\
&\leq4^{p-1}\left[E\|\xi\|^p+(c_1+c_2+c_3)T+(c_1+c_2+c_3)\int_0^tE\left(\sup\limits_{-\tau\leq
r\leq s}|y(r\wedge \theta_R)^-|^p\right)ds\right].
\end{align*}
Note that
\begin{equation*}
E\left(\sup\limits_{-\tau\leq s\leq t}| y(s\wedge \theta_R)^-|^p\right)\leq
E\|\xi\|^p+E\left(\sup\limits_{0\leq s\leq t}| y(s\wedge \theta_R)^-|^p\right).
\end{equation*}
Applying the Gronwall inequality and letting $R\to \infty,$  we then obtain
$$
E\left(\sup\limits_{-\tau\leq t\leq T}| y(t^-)|^p\right)\leq H(p).
$$
Since $T$ is any fixed positive number, the required assertion follows.
\end{proof}

In order to obtain our main results, we need to estimate the $p$th
moment of $y(s+\theta)-\bar{y}_s(\theta).$

\begin{lemma}\label{le3}
Let the conditions (\ref{eq5}) and (\ref{eq6}) hold. Then, for
$p\geq2$ and $s\in[0,T]$,
\begin{equation}\label{eq17}
E|y(s+\theta)-\bar{y}_s(\theta)|^p\leq \gamma\bigtriangleup,\ \ \ \
\ -\tau\leq\theta\leq0,
\end{equation}
where $\gamma$ is a positive constant which is independent of
$\bigtriangleup$.
\end{lemma}

\begin{proof}
Fix $s\in[0,T]$ and $\theta\in[-\tau,0]$. Let $k_s\in\{0, 1, 2,
\cdots, M-1\}$, $k_{\theta}\in\{-N, -N+1, \cdots, -1\}$ be the
integers for which $s\in[k_s\bigtriangleup, (k_s+1)\bigtriangleup)$,
$\theta\in[k_{\theta}\bigtriangleup, (k_{\theta}+1)\bigtriangleup)$,
 respectively. For convenience,
we write $v=s+\theta$ and $k_v=k_s+k_{\theta}$. Clearly, $0\leq
s-k_s\bigtriangleup<\bigtriangleup$ and $0\leq
\theta-k_{\theta}\bigtriangleup\leq\bigtriangleup$, so
\begin{equation*}
0\leq v-k_v\bigtriangleup<2\bigtriangleup.
\end{equation*}
Recalling the definition of $\bar{y}_s, s\in[0,T]$, we then yield
from \eqref{eq8} that
\begin{equation*}
\bar{y}_s(\theta)=\bar{y}_{k_s\bigtriangleup}(\theta)=\bar{y}(k_v\bigtriangleup)+
\frac{\theta-k_{\theta}\bigtriangleup}{\bigtriangleup}[\bar{y}((k_v+1)\bigtriangleup)-\bar{y}(k_v\bigtriangleup)],
\end{equation*}
which  implies
\begin{equation}\label{eq18}
\begin{split}
E|y(s+\theta)-\bar{y}_s(\theta)|^p\leq2^{p-1}E|\bar{y}((k_v+1)\bigtriangleup)-\bar{y}(k_v\bigtriangleup)|^p
+2^{p-1}E|y(v)-\bar{y}(k_v\bigtriangleup)|^p.
\end{split}
\end{equation}
For $k_v\leq-1$, it thus follows from \eqref{eq6} that
\begin{equation}\label{cy1}
E|\bar{y}((k_v+1)\bigtriangleup)-\bar{y}(k_v\bigtriangleup)|^p
\leq\beta\bigtriangleup.
\end{equation}
Note that, for some $\bar{H}:=\bar{H}(m,p)$,
\begin{equation}\label{eq28}
E|B(t)|^p\leq\bar{H}t^{\frac{p}{2}},\ \ \ \ \  t\geq0
\end{equation}
and, by the characteristic functions argument, for
$\bigtriangleup\in(0,1)$,
\begin{equation}\label{cy0}
E|\bigtriangleup N_k|^p\leq C\bigtriangleup,
\end{equation}
where $C$ is a positive constant which is independent of
$\bigtriangleup$. For $k_v\geq0,$ using \eqref{eq7} and noting
$g(\bar{y}_{k_v\bigtriangleup})$ and $B_{k_v},$
$h(\bar{y}_{k_v\bigtriangleup})$ and $N_{k_v}$ are independent,
respectively, we compute
\begin{equation*}
\begin{split}
E&|\bar{y}((k_v+1)\bigtriangleup)-\bar{y}(k_v\bigtriangleup)|^p\\
&\leq
3^{p-1}\left[E|f(\bar{y}_{k_v\bigtriangleup})|^p\bigtriangleup^p+E|g(\bar{y}_{k_v\bigtriangleup})|^pE|\bigtriangleup
B_{k_v}|^p +E|h(\bar{y}_{k_v\bigtriangleup})|^pE|\bigtriangleup
N_{k_v}|^p\right].
\end{split}
\end{equation*}
Taking \eqref{eq5} into consideration and applying Lemma \ref{le1},
we then obtain that for $\bigtriangleup\in(0,1)$
\begin{equation}\label{eq30}
E|\bar{y}((k_v+1)\bigtriangleup)-\bar{y}(k_v\bigtriangleup)|^p \leq
3^{p-1}2^{\frac{p}{2}-1}K^{\frac{p}{2}}(1+H(p))(1+\bar{H}+C)\bigtriangleup.
\end{equation}
Hence, in \eqref{eq18}
\begin{equation}\label{eq19}
\begin{split}
E|y(s+\theta)-\bar{y}_s(\theta)|^p&\leq[2^{p-1}\beta+3^{p-1}2^{\frac{3p}{2}-2}K^{\frac{p}{2}}(1+H(p))(1+\bar{H}+C)\bigtriangleup\\
&+2^{p-1}E|y(v)-\bar{y}(k_v\bigtriangleup)|^p.
\end{split}
\end{equation}
In what follows, we divide the following five cases to estimate the
second term on the right-hand side of \eqref{eq19}.

\noindent {\it Case 1:} $k_v\geq0$ and $0\leq
v-k_v\bigtriangleup<\bigtriangleup$. By \eqref{eq9}
\begin{equation*}
\begin{split}
E&|y(v)-\bar{y}(k_v\bigtriangleup)|^p\\
&=E|f(\bar{y}_{k_v\bigtriangleup})(v-k_v\bigtriangleup)+g(\bar{y}_{k_v\bigtriangleup})(B(v)-B(k_v\bigtriangleup))
+h(\bar{y}_{k_v\bigtriangleup})(N(v)-N(k_v\bigtriangleup))|^p\\
&\leq3^{p-1}E|f(\bar{y}_{k_v\bigtriangleup})|^p(v-k_v\bigtriangleup)^p+3^{p-1}E|g(\bar{y}_{k_v\bigtriangleup})|^pE|B(v)-B(k_v\bigtriangleup)|^p\\
&+3^{p-1}E|h(\bar{y}_{k_v\bigtriangleup})|^pE|N(v)-N(k_v\bigtriangleup)|^p.
\end{split}
\end{equation*}
Then, in the same way as \eqref{eq30} was done, we have for
$\bigtriangleup\in(0,1)$
\begin{equation*}
E|y(v)-\bar{y}(k_v\bigtriangleup)|^p
\leq3^{p-1}2^{\frac{p}{2}-1}K^{\frac{p}{2}}(1+H(p))(1+\bar{H}+C)\bigtriangleup.
\end{equation*}
\noindent {\it Case 2:} $k_v\geq0$ and $\bigtriangleup\leq
v-k_v\bigtriangleup<2\bigtriangleup$. It then follows easily that
\begin{equation*}
\begin{split}
E|y(v)-\bar{y}(k_v\bigtriangleup)|^p&=E|y(v)-\bar{y}((k_v+1)\bigtriangleup)+\bar{y}((k_v+1)\bigtriangleup)-\bar{y}(k_v\bigtriangleup)|^p\\
&\leq2^{p-1}E|y(v)-\bar{y}((k_v+1)\bigtriangleup)|^p+2^{p-1}E|\bar{y}((k_v+1)\bigtriangleup)-\bar{y}(k_v\bigtriangleup)|^p.
\end{split}
\end{equation*}
This, together with \eqref{eq30} and Case 1, leads to
\begin{equation*}
\begin{split}
E|y(v)-\bar{y}(k_v\bigtriangleup)|^p\leq3^{p-1}2^{\frac{3p}{2}-1}K^{\frac{p}{2}}(1+H(p))(1+\bar{H}+C)\bigtriangleup.
\end{split}
\end{equation*}
\noindent {\it Case 3:} $k_v=-1$ and $0\leq
v-k_v\bigtriangleup\leq\bigtriangleup$. In this case,
$-\bigtriangleup\leq v\leq0$. We then have from \eqref{eq6} that
\begin{equation*}
E|y(v)-\bar{y}(k_v\bigtriangleup)|^p\leq\beta\bigtriangleup.
\end{equation*}
\noindent {\it Case 4:} $k_v=-1$ and $\bigtriangleup\leq
v-k_v\bigtriangleup<2\bigtriangleup$. In such case, $0\leq
v<\bigtriangleup$. Case 1 and Case 2 can be used to estimate the
term
\begin{equation*}
\begin{split}
E|y(v)-\bar{y}(k_v\bigtriangleup)|^p&\leq2^{p-1}E|y(v)-\xi(0)|^p+2^{p-1}E|\xi(0)-\bar{y}((k_v\bigtriangleup)|^p\\
&\leq[2^{p-1}\beta+3^{p-1}2^{\frac{3p}{2}-2}K^{\frac{p}{2}}(1+H(p))(1+\bar{H}+C)]\bigtriangleup.
\end{split}
\end{equation*}
\noindent {\it Case 5:} $k_v\leq-2$. In this case, $v<0$. So, by
\eqref{eq6}
\begin{equation*}
E|y(v)-\bar{y}(k_v\bigtriangleup)|^p\leq2\beta\bigtriangleup.
\end{equation*}
Combining case 1 to case 5,  we therefore complete the proof.
\end{proof}

The following Theorem will tell us the error of the $p$th moment
between the true solution and numerical solution under global
Lipschitz condition.

\begin{theorem}\label{le4}
Under the conditions \eqref{eq4} and \eqref{eq6}, for $p\geq2$,
\begin{equation}\label{cy2}
E\left(\sup\limits_{0\leq t\leq T}|
x(t)-y(t)|^p\right)\leq\delta_1L^{\frac{p}{2}}
e^{\delta_2L^{\frac{p}{2}}}\bigtriangleup,
\end{equation}
where $\delta_1, \delta_2$ are constants which are independent of $\bigtriangleup.$
\end{theorem}

\begin{proof}
It is easy to see from \eqref{eq1} and \eqref{eq9} that for any
$t_1\in[0,T]$
\begin{equation}\label{1}
\begin{split}
E\left(\sup\limits_{0\leq t\leq
t_1}|x(t)-y(t)|^p\right)&\leq3^{p-1}E\left(\sup\limits_{0\leq t\leq
t_1}\left|\int_0^tf(x_{s^-})-f(\bar{y}_{s^-})ds\right|^p\right)\\
&+3^{p-1}E\left(\sup\limits_{0\leq t\leq
t_1}\left|\int_0^tg(x_{s^-})-g(\bar{y}_{s^-})dB(s)\right|^p\right)\\
&+3^{p-1}E\left(\sup\limits_{0\leq t\leq
t_1}\left|\int_0^th(x_{s^-})-h(\bar{y}_{s^-})dN(s)\right|^p\right)\\
&:=I_1+I_2+I_3.
\end{split}
\end{equation}
In the sequel, we estimate these terms respectively. By the H\"older
inequality, \eqref{eq4} and Lemma \ref{le3},
\begin{equation*}
\begin{split}
I_1&\leq3^{p-1}T^{p-1}\int_0^{t_1}E|f(x_s)-f(\bar{y}_s)|^pds\\
&\leq6^{p-1}T^{p-1}\int_0^{t_1}E|f(x_s)-f(y_s)|^pds+6^{p-1}T^{p-1}\int_0^{t_1}E|f(y_s)-f(\bar{y}_s)|^pds\\
&\leq6^{p-1}T^{p-1}\int_0^{t_1}E\left(\int_{-\tau}^0|x(s+\theta)-y(s+\theta)|^2d\mu(\theta)\right)^{\frac{p}{2}}ds\\
&+6^{p-1}T^{p-1}\int_0^{t_1}E\left(\int_{-\tau}^0|y(s+\theta)-\bar{y}_s(\theta)|^2d\mu(\theta)\right)^{\frac{p}{2}}ds\\
&\leq6^{p-1}T^{p-1}L^{\frac{p}{2}}\int_0^{t_1}E\left(\sup\limits_{0\leq
r\leq s}|x(r)-y(r)|^p\right)ds
+6^{p-1}T^{p-1}L^{\frac{p-2}{2}}\int_0^{t_1}\int_{-\tau}^0E|y(s+\theta)-\bar{y}_s(\theta)|^pd\mu(\theta)ds\\
&\leq6^{p-1}T^pL^{\frac{p}{2}}\gamma\bigtriangleup+6^{p-1}T^{p-1}L^{\frac{p}{2}}\int_0^{t_1}E\left(\sup\limits_{0\leq
r\leq s}|x(r)-y(r)|^p\right)ds.
\end{split}
\end{equation*}
Now, the Burkholder-Davis-Gundy inequality \cite[Theorem 7.3,
p40]{m97}, \eqref{eq4} and Lemma \ref{le3} also give that, for some
positive constant $C_p$,
\begin{equation}\label{2}
\begin{split}
I_2&\leq 3^{p-1}C_pE\left(\int_0^{t_1}|g(x_s)-g(\bar{y}_s)|^2ds\right)^{\frac{p}{2}}\\
&\leq3^{p-1}T^{\frac{p-2}{2}}C_p\int_0^{t_1}E|g(x_s)-g(\bar{y}_s)|^pds\\
&\leq6^{p-1}T^{\frac{p-2}{2}}C_p\left[L^{\frac{p}{2}}\int_0^{t_1}E\left(\sup\limits_{0\leq
r\leq
s}|x(r)-y(r)|^p\right)ds+L^{\frac{p-2}{2}}\int_0^{t_1}\int_{-\tau}^0E|y(s+\theta)-\bar{y}_s(\theta)|^pd\mu(\theta)ds\right]\\
&\leq6^{p-1}\gamma
T^{\frac{p}{2}}C_pL^{\frac{p}{2}}\bigtriangleup+6^{p-1}T^{\frac{p-2}{2}}C_pL^{\frac{p}{2}}\int_0^{t_1}E\left(\sup\limits_{0\leq
r\leq s}|x(r)-y(r)|^p\right)ds.
\end{split}
\end{equation}
In the same way as \eqref{2} was done, together with the
Burkholder-Davis-Gundy inequality \cite[Theorem 48, p193]{p04}, we
can deduce from \eqref{eq4} that, for some positive constant
$\bar{C}_p$,
\begin{equation*}
\begin{split}
I_3&\leq6^{p-1}E\left(\sup\limits_{0\leq t\leq
t_1}\left|\int_0^th(x_{s^-})-h(\bar{y}_{s^-})d\tilde{N}(s)\right|^p+\lambda^p\sup\limits_{0\leq
t\leq t_1}\left|\int_0^th(x_s)-h(\bar{y}_s)ds\right|^p\right)\\
&\leq6^{p-1}(\bar{C}_pT^{\frac{p-2}{2}}\lambda^{\frac{p}{2}}+\lambda^pT^{p-1})\int_0^{t_1}E|h(x_s)-h(\bar{y}_s)|^pds\\
&\leq12^{p-1}(\bar{C}_pT^{\frac{p-2}{2}}\lambda^{\frac{p}{2}}+\lambda^pT^{p-1})L^{\frac{p}{2}}\int_0^{t_1}E\left(\sup\limits_{0\leq
r\leq
s}|x(r)-y(r)|^p\right)ds\\
&+12^{p-1}(\bar{C}_pT^{\frac{p-2}{2}}\lambda^{\frac{p}{2}}+\lambda^pT^{p-1})L^{\frac{p-2}{2}}\int_0^{t_1}\int_{-\tau}^0E|y(s+\theta)-\bar{y}_s(\theta)|^pd\mu(\theta)ds\\
&\leq12^{p-1}\gamma(\bar{C}_pT^{\frac{p}{2}}\lambda^{\frac{p}{2}}+\lambda^pT^p)L^{\frac{p}{2}}\bigtriangleup+12^{p-1}(\bar{C}_pT^{\frac{p-2}{2}}\lambda^{\frac{p}{2}}+\lambda^pT^{p-1})L^{\frac{p}{2}}\int_0^{t_1}E\left(\sup\limits_{0\leq
r\leq s}|x(r)-y(r)|^p\right)ds.
\end{split}
\end{equation*}
Therefore
\begin{equation*}
E\left(\sup\limits_{0\leq t\leq
t_1}|x(t)-y(t)|^p\right)\leq\delta_1L^{\frac{p}{2}}\bigtriangleup
+\delta_2L^{\frac{p}{2}}\int_0^{t_1}E\left(\sup\limits_{0\leq r\leq
s}|x(r)-y(r)|^p\right)ds,
\end{equation*}
where $\delta_1=6^{p-1}\gamma
T^{\frac{p}{2}}\Big(T^{\frac{p}{2}}+C_p+2^{p-1}\bar{C}_p\lambda^{\frac{p}{2}}+\lambda^pT^{\frac{p}{2}}\Big)$
and
$\delta_2=6T^{\frac{p-2}{2}}\Big(T^{\frac{p}{2}}+C_p+2^{p-1}\bar{C}_p\lambda^{\frac{p}{2}}+\lambda^pT^{\frac{p}{2}}\Big)$.
The desired assertion thus follows from the Gronwall inequality.
\end{proof}

\begin{remark}
The result of  Theorem \ref{le4} tells us
\begin{equation}\label{cy3}
E\left(\sup\limits_{0\leq t\leq T}|
x(t)-y(t)|^2\right)\leq \delta_3 L
e^{\delta_4L}\bigtriangleup,
\end{equation}
where $\delta_3, \delta_4$ are constants which are independent of $\bigtriangleup$
under the global Lipschitz condition \eqref{eq4}.  This means that the order  of the
mean-square convergence is $1/2$,
 while Eq. \eqref{cy2} tells us that the order of the $p$th moment convergence is $1/p\  (p\ge 2)$.  In other words, the lower moment has a  better convergence
 rate for the SFDEs with jumps, whence it is best in practice to use the mean-square convergence. This is significantly different from the result on SFDEs without jumps.
 Letting $h\equiv 0$ in \eqref{eq1},
 i.e. there is no jumps, we have already known that for $p\ge 2$ (see \cite{ym08})
\begin{equation*}
E\left(\sup\limits_{0\leq t\leq T}|
x(t)-y(t)|^p\right)\leq \hat C_1 \bigtriangleup^{p/2},
\end{equation*}
where $\hat C_1$ is a  constant independent of
$\bigtriangleup.$   This means that the order of the $p$th moment convergence is $1/2$
 for all $p\ge 2$.
 Why is there a significant
difference? Actually, it is due to the following fact: all moments of the Poisson increments
$N((k+1)\bigtriangleup)-N(k\bigtriangleup)$ have the same order of
$\bigtriangleup$ $( see (\ref{cy0}) )$,
 while the moments of increments $\bigtriangleup B_k=B((k+1)\bigtriangleup)-B(k\bigtriangleup)$
 have different orders, namely
  $E |\bigtriangleup B_k|^{2n}= O(\bigtriangleup^n)$ and
$E |\bigtriangleup B_k|^{2n+1}= 0.$
\end{remark}

\section{Rate of convergence under local Lipschitz condition}

In this section, we shall discuss the rate of convergence of EM
numerical solutions to \eqref{eq1} under the following {\it local
Lipschitz condition}.

\begin{enumerate}
\item[\textmd{(\bf H3)}] (Local Lipschitz condition) For each integer $j\geq1$, there is a left-continuous nondecreasing function
$\mu_j:[-\tau,0]\rightarrow R_+$ such that
\begin{equation}\label{eq20}
|f(\varphi)-f(\psi)|^2\vee|g(\varphi)-g(\psi)|^2\vee|h(\varphi)-h(\psi)|^2\leq
\int_{-\tau}^0|\varphi(\theta)-\psi(\theta)|^2d\mu_j(\theta),
\end{equation}
\end{enumerate}
for those $\varphi,\psi \in \hat D([-\tau,0];R^n)$ with
$\|\varphi\|\vee\|\psi\|\leq j$.

\begin{enumerate}
\item[\textmd{(\bf H4)}] (Linear growth condition) Assume that there
is a constant $h>0$ such that, for $\varphi\in \hat D([-\tau,0];R^n)$,
\begin{equation}\label{eq31}
|f(\varphi)|^2\vee |g(\varphi)|^2\vee |h(\varphi)|^2\leq
h(1+\|\varphi\|^2).
\end{equation}
\end{enumerate}

\begin{remark}
Under the conditions \eqref{eq20} and \eqref{eq31}, for any initial
data $\xi\in D^{b}_{{\mathcal F}_{0}}([-\tau,0];R^n)$, \eqref{eq1}
admits a unique solution $x(t), t\in[0,T]$ by using  the standard truncation procedure (see \cite[Theorem 3.4,
p56]{m97} ). Moreover,  \eqref{eq20} implies for those $\varphi,\psi \in \hat D([-\tau,0];R^n)$ with
$\|\varphi\|\vee\|\psi\|\leq j$
\begin{equation}
|f(\varphi)-f(\psi)|^2\vee|g(\varphi)-g(\psi)|^2\vee|h(\varphi)-h(\psi)|^2\leq
L_j \|\varphi-\psi\|^2,
\end{equation}
where $L_j=\mu_j(0)-\mu_j(-\tau)$.
\end{remark}

\begin{theorem}
Let
conditions \eqref{eq6}, \eqref{eq20} and \eqref{eq31} hold. If there exist  positive
constant $\alpha$  and  $\tilde \varepsilon \in (0, 1)$ such that the local
Lipschitz constant obeys
\begin{equation}\label{eq21}
L_j^{1+\tilde \varepsilon}\leq\alpha\log j,
\end{equation}
then
\begin{equation}\label{eq22}
E\left(\sup\limits_{0\leq t\leq T}| x(t)-y(t)|^2\right)=
O(\bigtriangleup^{\frac{2}{2+\epsilon}}),
\end{equation}
where $\epsilon\in (0, \tilde \varepsilon)$ is an arbitrarily fixed
small positive number.
\end{theorem}

\begin{proof}
Let $j\geq1$ be an integer, and let $S_j=\{x\in R^n: |x|\leq j\}$.
Define the projection $\pi_j:R^n\rightarrow S_j$ by
\begin{equation*}
\pi_j(x)=\frac{j\wedge |x|}{|x|}x,
\end{equation*}
where we set $\pi_j(0)=0$ as usual. It is easy to see that for all
$x,y\in R^n$
\begin{equation*}
|\pi_j(x)-\pi_j(y)|\leq|x-y|.
\end{equation*}
Define the operator $\bar{\pi}_j: \hat D([-\tau,0];R^n)\rightarrow
\hat D([-\tau,0];R^n)$ by
\begin{equation*}
\bar{\pi}_j(\varphi)=\{\pi_j(\varphi(\theta)):-\tau\leq\theta\leq0\}.
\end{equation*}
Clearly,
\begin{equation*}
\|\bar{\pi}_j(\varphi)\|\leq j,\ \ \ \ \ \ \ \ \forall \varphi\in
\hat D([-\tau,0];R^n).
\end{equation*}
Define the truncation functions $f_j:\hat D([-\tau,0];R^n)\rightarrow
R^n$, $g_j:\hat D([-\tau,0];R^n)\rightarrow R^{n\times m}$ and
$h_j:  \hat D([-\tau,0];R^n)\rightarrow R^n$ by
\begin{equation}\label{eq23}
f_j(\varphi)=f(\bar{\pi}_j(\varphi)),\ \ \ \ \
g_j(\varphi)=g(\bar{\pi}_j(\varphi)),\ \ \ \ \
h_j(\varphi)=h(\bar{\pi}_j(\varphi)),
\end{equation}
respectively. Then, by \eqref{eq20}, for any $\varphi,\psi\in
\hat D([-\tau,0];R^n)$,
\begin{equation}
\begin{split}
&|f_j(\varphi)-f_j(\psi)|^2\vee|g_j(\varphi)-g_j(\psi)|^2\vee|h_j(\varphi)-h_j(\psi)|^2\\
&\leq|f(\bar{\pi}_j(\varphi))-f(\bar{\pi}_j(\psi))|^2\vee|g(\bar{\pi}_j(\varphi))-g(\bar{\pi}_j(\psi))|^2
\vee|h(\bar{\pi}_j(\varphi))-h(\bar{\pi}_j(\psi))|^2\\
&\leq\int_{-\tau}^0|\pi_j(\varphi(\theta))-\pi_j(\psi(\theta))|^2d\mu_j(\theta)\\
&\leq\int_{-\tau}^0|\varphi(\theta)-\psi(\theta)|^2d\mu_j(\theta).
\end{split}
\end{equation}
That is, $f_j$, $g_j$ and $h_j$ satisfy the global Lipschitz
condition.  For $t\in[0,T]$, let $x^j(t)$ be the solution to the
following SFDE with jumps
\begin{equation*}
dx^j(t)=f_j(x^j_{t^-})dt+g_j(x^j_{t^-})dB(t)+h_j(x^j_{t^-})dN(t)
\end{equation*}
with the initial data $x^j_0=\xi$ and $y^j(t)$ be the corresponding
continuous-time EM solution with the stepsize $\bigtriangleup$. By
Theorem \ref{le4} for any sufficiently small$\epsilon \in (0,\tilde
\varepsilon)$
\begin{equation*}
E\left(\sup\limits_{0\leq t\leq T}|
x^j(t)-y^j(t)|^{2+\epsilon}\right)\leq\delta_1L_j^{1+\epsilon/2}e^{\delta_2L_j^{1+\epsilon/2}}\bigtriangleup.
\end{equation*}
Furthermore, by \eqref{eq21} (here we assume $L_j\ge 1$ without any
loss of generality),
\begin{equation}\label{eq24}
E\left(\sup\limits_{0\leq t\leq T}| x^j(t)-y^j(t)|^{2+\epsilon}\right)\leq
e^{(\delta_1+\delta_2)L_j^{1+\epsilon/2}}\bigtriangleup\leq
j^{\alpha(\delta_1+\delta_2)}\bigtriangleup.
\end{equation}
Set
\begin{equation*}
\hat{x}(T)=\sup\limits_{0\leq t\leq T}| x(t)|\ \ \ and \ \ \
\hat{y}(T)=\sup\limits_{0\leq t\leq T}| y(t)|.
\end{equation*}
For any integer $j\geq1$, define stopping time
\begin{equation*}
\tau_j=T\wedge\inf\{t\in[0,T]: \| x^j_t\|\vee\| y^j_t\|\ge j\}.
\end{equation*}
It is easy to see that $\| x^j_{s^-}\|\leq j$ for any $0\leq
s<\tau_j$. Then, combining \eqref{eq23} gives that for any $0\leq
s<\tau_j$
\begin{equation*}
f_j(x^j_{s^-})=f\left(\frac{\| x^j_{s^-}\|\wedge j}{\|
x^j_{s^-}\|}x^j_{s^-}\right)=f\left(\frac{\| x^j_{s^-}\|\wedge
(j+1)}{\|
x^j_{s^-}\|}x^j_{s^-}\right)=f_{j+1}(x^j_{s^-})=f(x^j_{s^-}).
\end{equation*}
Similarly,
\begin{equation*}
g_j(x^j_{s^-})=g_{j+1}(x^j_{s^-})=g(x^j_{s^-}),\ \ \ \
h_j(x^j_{s^-})=h_{j+1}(x^j_{s^-})=h(x^j_{s^-}).
\end{equation*}
While on $0\leq t<\tau_j$
\begin{equation*}
\begin{split}
x^j(t)&=\xi(0)+\int_0^tf_j(x^j_{s^-})ds+\int_0^tg_j(x^j_{s^-})dB(s)+\int_0^th_j(x^j_{s^-})dN(s)\\
&=\xi(0)+\int_0^tf_{j+1}(x^j_{s^-})ds+\int_0^tg_{j+1}(x^j_{s^-})dB(s)+\int_0^th_{j+1}(x^j_{s^-})dN(s)\\
&=\xi(0)+\int_0^tf(x^j_{s^-})ds+\int_0^tg(x^j_{s^-})dB(s)+\int_0^th(x^j_{s^-})dN(s).
\end{split}
\end{equation*}
Consequently, we must have that
\begin{equation*}
x(t)=x^j(t)=x^{j+1}(t)
\end{equation*}
on $0\leq t<\tau_j$. Likewise, we can also derive that
\begin{equation*}
y(t)=y^j(t)=y^{j+1}(t)
\end{equation*}
for $0\leq t<\tau_j$. These imply that $\tau_j$ is non-decreasing
and, by Lemma \ref{le1},
$\lim\limits_{j\rightarrow\infty}\tau_j=T$ a.s. Let $\tau_0=0$ and
compute, for $t\in[0,T]$,
\begin{equation*}
\begin{split}
| x(t)-y(t)|^2&=\sum\limits_{j=1}^{\infty}|
x(t)-y(t)|^2I_{[\tau_{j-1}\leq t<\tau_j)}\\
&=\sum\limits_{j=1}^{\infty}|
x^j(t)-y^j(t)|^2I_{[\tau_{j-1}\leq t<\tau_j)}\\
&\leq\sum\limits_{j=1}^{\infty}| x^j(t)-y^j(t)|^2I_{[j-1\leq
\hat{x}(T)\vee\hat{y}(T)]}.
\end{split}
\end{equation*}
Therefore, by the H\"older inequality
\begin{equation}\label{eq25}
\begin{split}
&E\left(\sup\limits_{0\leq t\leq T}|
x(t)-y(t)|^2\right)\\
&\leq\sum\limits_{j=1}^{\infty}\left(E\left(\sup\limits_{0\leq
t\leq T}| x^j(t)-y^j(t)|^{2+\epsilon}\right)\right)^{\frac{2}{2+\epsilon}}\Big(EI_{[j-1\leq
\hat{x}(T)\vee\hat{y}(T)]}\Big)^{\frac{\epsilon}{2+\epsilon}}\\
&\leq\sum\limits_{j=1}^{\infty}\left(E\left(\sup\limits_{0\leq t\leq
T}| x^j(t)-y^j(t)|^{2+\epsilon}\right)\right)^{\frac{2}{2+\epsilon}}[P(j-1\leq
\hat{x}(T)\vee\hat{y}(T))]^{\frac{\epsilon}{2+\epsilon}}.
\end{split}
\end{equation}
On the other hand, for any $q\geq2$, we obtain from Lemma \ref{le1}
\begin{equation}  \label{eq26}
P(j-1\leq \hat{x}(T)\vee\hat{y}(T))\leq\frac{E| \hat{x}(T)|^q+E|
\hat{y}(T)|^q}{(\frac{j}{2})^q}\leq\frac{2H(q)}{(\frac{j}{2})^q}
\end{equation}
with $j\geq2$. Substituting \eqref{eq24} and \eqref{eq26} into
\eqref{eq25}, one has
\begin{equation}\label{eq27}
E\left(\sup\limits_{0\leq t\leq T}| x(t)-y(t)|^2\right)
\leq\left(1+2^{\frac{q\epsilon}{2+\epsilon}}(2H(q))^{\frac{\epsilon}{2+\epsilon}}\sum\limits_{j=2}^{\infty}
j^{\frac{2\alpha(\delta_1+\delta_2)-q\epsilon}{2+\epsilon}}\right)\bigtriangleup^{\frac{2}{2+\epsilon}}.
\end{equation}
For any fixed $\epsilon>0$ letting $q$ be sufficiently large for
$$
q\geq\frac{\alpha(\delta_1+\delta_2)+2(2+\epsilon)}{\epsilon},
$$ we see that the right-hand side
of \eqref{eq27} is convergent, whence the desired assertion
\eqref{eq22} follows.
\end{proof}

\begin{remark}
Under the local Lipschitz condition, Mao \cite{m03} showed the
strong convergence of the numerical solutions to SFDEs without
jumps, and the rate of convergence was revealed under the global
Lipschitz condition.  In the present paper, under the local
Lipschitz condition, we reveal the rate of convergence for the
numerical solutions to SFDEs with jumps.  The rate of convergence
for jump processes \eqref{eq1} we revealed here is $1/(2+\epsilon)$
(closed to $1/2$) under the logarithm growth condition \eqref{eq21}.
This is different from the rate of convergence for the diffusion
processes (without jumps) which was studied in \cite{ym08}, where it
was  shown that the rate of convergence is still $1/2$  under the
logarithm growth condition.  The reason for such a difference has
already been pointed out in Remark 3.1.
\end{remark}

\section{Appendix: an existence-and-uniqueness theorem}
To make our paper self-contained, in this section we shall discuss
the existence and uniqueness of solutions to \eqref{eq1} under the
assumption $(H1)$.

\begin{theorem}
Under the conditions \eqref{eq2}, there exists a unique solution
$x(t), t\in[0,T],$ to \eqref{eq1} for any initial data $\xi\in
D^{b}_{{\mathcal F}_{0}}([-\tau,0];R^n)$.
\end{theorem}

\begin{proof}
Since our proof is an application of the proof for the case without
jumps in \cite[Theorem 2.2, p150]{m97}, we here give only a sketch
for the proof of jump case.

\noindent{\it Uniqueness.} Let $x(t)$ and $\bar{x}(t)$ be two
solutions to \eqref{eq1} on $[0,T]$. Noting from \eqref{eq1} that
\begin{equation*}
x(t)-\bar{x}(t)=\int_0^t[f(x_{s^-})-f(\bar{x}_{s^-})]ds+\int_0^t[g(x_{s^-})-g(\bar{x}_{s^-})]dB(s)+\int_0^t[h(x_{s^-})-h(\bar{x}_{s^-})]dN(s)
\end{equation*}
and $\tilde{N}(t)=N(t)-\lambda t$ is a martingale measure for
$t\in[0,T]$, along with \eqref{eq4} we have
\begin{equation*}
E\left(\sup\limits_{0\leq t\leq
T}|x(t)-\bar{x}(t)|^2\right)\leq3L(T+4+8\lambda+2\lambda^2T)\int_0^TE\left(\sup\limits_{0\leq
r\leq s}|x(r)-\bar{x}(r)|^2\right)ds.
\end{equation*}
By the Gronwall inequality
\begin{equation*}
E\left(\sup\limits_{0\leq t\leq T}|x(t)-\bar{x}(t)|^2\right)=0,
\end{equation*}
which implies that $x(t)=\bar{x}(t)$ for $t\in[0,T]$ almost surely.
The uniqueness has been proved.

\noindent{\it Existence.} Define $x^0_0=\xi$ and $x^0(t)=\xi(0)$ for
$0\leq t\leq T$. For each $n=1,2,\cdots,$ set $x_0^n=\xi$ and
define, by the Picard iterations,
\begin{equation}\label{eq32}
x^n(t)=\xi(0)+\int_0^tf(x_{s^-}^{n-1})ds+\int_0^tg(x_{s^-}^{n-1})dB(s)+\int_0^th(x_{s^-}^{n-1})dN(s)
\end{equation}
for $t\in[0,T]$. It also follows from \eqref{eq32} that for any
integer $k\geq1$
\begin{equation*}
E\left(\sup\limits_{0\leq t\leq T}|x^n(t)|^2\right)\leq
\bar{c}_1+\bar{c}_2\int_0^TE\left(\sup\limits_{0\leq r\leq
s}|x^{n-1}(r)|^2\right)ds,
\end{equation*}
where
$\bar{c}_1=4[E\|\xi\|^2+K(T^2+4T)+(8\lambda+2\lambda^2T)T]+\bar{c}_2TE\|\xi\|^2$
and $\bar{c}_2=4K[T+4+8\lambda+2\lambda^2T]$. This further implies
that
\begin{equation*}
\begin{split}
\max\limits_{1\leq n\leq k}E\left(\sup\limits_{0\leq t\leq
T}|x^n(t)|^2\right)&\leq
\bar{c}_1+\bar{c}_2\int_0^T\max\limits_{1\leq n\leq
k}E\left(\sup\limits_{0\leq s\leq t}|x^{n-1}(s)|^2\right)dt.
\end{split}
\end{equation*}
Observing
\begin{equation*}
\begin{split}
&\max\limits_{1\leq n\leq k}E\left(\sup\limits_{0\leq s\leq
t}|x^{n-1}(s)|^2\right)\\
&=\max\left\{E|\xi(0)|^2,E\left(\sup\limits_{0\leq s\leq
t}|x^1(s)|^2\right),\cdots, E\left(\sup\limits_{0\leq s\leq
t}|x^{k-1}(s)|^2\right) \right\}\\
&\leq\max\left\{E\|\xi\|^2,E\left(\sup\limits_{0\leq s\leq
t}|x^1(s)|^2\right),\cdots, E\left(\sup\limits_{0\leq s\leq
t}|x^{k-1}(s)|^2\right),E\left(\sup\limits_{0\leq s\leq
t}|x^k(s)|^2\right) \right\}\\
&\leq E\|\xi\|^2+\max\limits_{1\leq n\leq
k}E\left(\sup\limits_{0\leq s\leq t}|x^n(s)|^2\right),
\end{split}
\end{equation*}
we hence deduce that
\begin{equation*}
\begin{split}
\max\limits_{1\leq n\leq k}E\left(\sup\limits_{0\leq t\leq
T}|x^n(t)|^2\right)
\leq\bar{c}_1+\bar{c}_2TE\|\xi\|^2+\bar{c}_2\int_0^T\max\limits_{1\leq
n\leq k}E\left(\sup\limits_{0\leq s\leq t}|x^n(s)|^2\right)dt.
\end{split}
\end{equation*}
The Gronwall inequality implies
\begin{equation*}
\max\limits_{1\leq n\leq k}E\left(\sup\limits_{0\leq t\leq
T}|x^n(t)|^2\right)\leq(\bar{c}_1+\bar{c}_2TE\|\xi\|^2)e^{\bar{c}_2T}.
\end{equation*}
Since $k$ is arbitrary, we must have for $n\geq1$
\begin{equation*}
E\left(\sup\limits_{0\leq t\leq
T}|x^n(t)|^2\right)\leq(\bar{c}_1+\bar{c}_2TE\|\xi\|^2)e^{\bar{c}_2T}.
\end{equation*}
Next, by \eqref{eq32}
\begin{equation}\label{eq33}
\begin{split}
E\left(\sup\limits_{0\leq t\leq
T}|x^1(t)-x^0(t)|^2\right)&\leq3K(T+4+8\lambda+2\lambda^2T)\int_0^T(1+E\|x^0_s\|^2)ds\\
&\leq3KT(T+4+8\lambda+2\lambda^2T)(1+E\|\xi\|^2):=\bar{C}.
\end{split}
\end{equation}
We now claim that for $n\geq0$
\begin{equation}\label{eq34}
E\left(\sup\limits_{0\leq s\leq
t}|x^{n+1}(s)-x^n(s)|^2\right)\leq\frac{\bar{C}M^nt^n}{n!},\ \ \ \
\ 0\leq t\leq T,
\end{equation}
where $M=3K(T+4+8\lambda+2\lambda^2T)$. We shall show this by
induction. In view of \eqref{eq33} we see that \eqref{eq34} holds
whenever $n=0$. Now, assume that \eqref{eq34} holds for some
$n\geq0$. Then,
\begin{equation*}
\begin{split}
E\left(\sup\limits_{0\leq s\leq
t}|x^{n+2}(s)-x^{n+1}(s)|^2\right)&\leq
M\int_0^tE\|x^{n+1}_s-x^n_s\|^2ds\\
&\leq M\int_0^tE\left(\sup\limits_{0\leq r\leq
s}|x^{n+1}(r)-x^n(r)|^2\right)ds\\
&\leq
M\int_0^t\frac{\bar{C}M^ns^n}{n!}ds=\frac{\bar{C}M^{n+1}t^{n+1}}{(n+1)!}.
\end{split}
\end{equation*}
Following the
proof of \cite[Theorem 3.1, p55]{m97}, we can show that for almost all $\omega \in \Omega$ there exists a positive integer $n_0=n_0(\omega)$ such that
\begin{align}
\sup_{0\le s\le T}|x^{n+1}(s)-x^{n}(s)|\le \frac{1}{2^n}\  \mbox{
whenever } n\ge n_0(\omega).
\end{align}
This implies $\{x^{n}(\cdot)\}_{n\ge 1}$ is a Cauchy sequence under
$\sup|\cdot|.$ However, since our space $D([0, T]; R^n)$ is not a
complete space under $\sup|\cdot|,$ we do not know whether
$\{x^{n}(\cdot)\}_{n\ge 1}$ has a limit in $D([0, T]; R^n).$ In order
for $D([0,T]; R^n)$ is complete, we need to define the
following metric (see \cite[Chapter 3]{bil68}). Let $\Lambda$ denote
the class of strictly increasing, continuous mapping of $[0, T]$
onto itself and
$$
\Lambda_{\epsilon}^*=\left\{\lambda\in \Lambda: \sup_{s\neq t}\left|\log\frac{\lambda(t)-\lambda(s)}{t-s}\right|\le \epsilon\right\},
$$
define
\begin{align*}
d(\xi, \zeta)=\inf\{\epsilon>0: \exists \lambda\in \Lambda_{\epsilon}^* \mbox{ such that }
\sup_{t\in [0, T]}|\xi(t)-\zeta(\lambda(t))|\le \epsilon\}.
\end{align*}
$d(\cdot, \cdot)$ is called the Skorohod metric, by \cite[Theorem 14.2, p115]{bil68}  we know that  $D([0, T]; R^n)$ is complete in the metric $d.$ Taking $\lambda(t)=t,$ we can see $\{x^{n}(\cdot)\}_{n\ge 1}$ is a Cauchy sequence under $d.$ Therefore there exists unique
$x(t),
t\in[0,T] \in D([0, T]; R^n)$ such that $d(x^{n}(\cdot), x(\cdot))\rightarrow 0$ as $n\rightarrow \infty.$ Taking the limit in \eqref{eq32}, we then can show that $x(t)$ is  the  solution of \eqref{eq1}.
\end{proof}

\end{document}